\documentclass[conference]{IEEEtran}

\usepackage{algorithm}
\usepackage{algorithmic}
\usepackage[cmex10]{amsmath}
\usepackage{amsfonts, amssymb}
\usepackage{bm} 
\usepackage{enumitem}
\usepackage[caption=false,font=footnotesize]{subfig}
\usepackage[pdftex]{graphicx}
\usepackage{mathrsfs}
\usepackage[dvipsnames]{xcolor}
\usepackage{mathtools}
\usepackage[noadjust]{cite}

\graphicspath{{./figs/}}
\DeclareGraphicsExtensions{.png,.jpg,.pdf,.eps}

\newcommand{\norm}[1]{\left\lVert#1\right\rVert} 
\newcommand{\real}{\mathbb{R}} 
\newcommand{\ones}{\mathbf{1}}

\newcommand{\trans}{^\top}
\DeclareMathOperator*{\argmin}{arg\,min} 

\begin{document}

\title{Adaptive Coordination of Distributed Energy Resources in Lossy Power Distribution Systems}

\author{\IEEEauthorblockN{Hanchen Xu, Alejandro~D.~Dom\'{i}nguez-Garc\'{i}a, and Peter~W.~Sauer}
\IEEEauthorblockA{Department of Electrical and Computer Engineering\\
University of Illinois at Urbana-Champaign\\
Urbana, IL 61801\\
Email: \{hxu45, aledan, psauer\}@illinois.edu}
}

\maketitle

\begin{abstract}
This paper is concerned with the problem of coordinating a set of distributed energy resources (DERs) in a lossy power distribution system to provide frequency regulation services to a bulk power grid with the explicit consideration of system losses.
To this end, we formulate the problem as an optimization problem, the objective of which is to minimize some cost function subject to a set of constraints.
The formulation requires knowledge of incremental total system losses, which we approximate using the so-called loss factors (LFs) that explicitly capture the impacts of both active and reactive power injections on system losses.
The LFs are estimated recursively using power injection measurements; thus, they are adaptive to various phenomena that impact the power system operation. 
Numerical simulation on a 33-bus distribution test feeder validated the effectiveness of the proposed framework.
\end{abstract}



\section{Introduction} \label{sec:intro}

Frequency regulation services are by-and-large provided by conventional synchronous generators.
However, with deepening penetration of renewable-based generation resources, synchronous generators alone may be insufficient to meet regulation requirements \cite{wind_impacts}.
Moreover, conventional synchronous generators may suffer from poor performances when the regulation signal changes fast.
To overcome the aforementioned challenges, distributed energy resources (DERs), such as energy storage resources, are allowed to provide such services \cite{ES}.
In this context in the US, driven by  Order 755 from the Federal Energy Regulatory Commission, performance-based regulation markets have emerged as an effective means to incentivize the provision of high-quality frequency regulation services from resources including conventional generators and DERs \cite{performance}.
Resources are incentivized to track the instructed regulation signal accurately in a performance-based regulation market since otherwise they will incur loss of payments \cite{performance}. 

Compared to conventional generators, DERs typically enjoy much faster-responding capabilities, and can potentially track the regulation signal better.
Yet, the individual frequency regulation capability of the DER is usually small \cite{DER_coord2}.
As such, appropriate coordination of DERs is required in order for them to collectively provide some amount of regulation power in the form of an incremental change (with respect to some nominal value) in the power exchanged by the distribution feeder with the rest of the system; this incremental change is specified by the regulation signal sent by the bulk system operator.
The problem of optimally allocating the regulation power among the set of DERs subject to capacity limits is referred to as the optimal DER coordination problem (ODCP).

The general ODCP that concerns with the allocation of resources (active or reactive power) has been widely studied \cite{DER_coord2, DER_model1}. 
For example, a distributed coordination scheme is proposed for DERs to provide ancillary services in lossless systems in \cite{DER_coord2}.
Authors in \cite{DER_model1} present a framework to coordinate the DERs to provide the primary frequency control using a  droop-control scheme, based on a power flow model with explicit consideration of losses.
However, neglecting system losses will result in a solution that fail to provide the exact amount regulation power at the power distribution feeder.
Also, the dependence on model parameters obtained offline and thus lack adaptivity to changes in system conditions as evidenced in \cite{meas_est, meas_est2}.
Moreover, in situations where the model parameters are unavailable, these approaches simply cannot work.

In this paper, we address the ODCP by taking an alternative approach that explicitly takes into account system losses yet without reliance on system models.
To this end, we formulate the ODCP as an optimization problem that aims to minimize some cost function and is constrained by the power balance equation and DER capacity limits.
In order to provide the exact amount of regulation power, it is necessary to explicitly consider the impacts of system losses, which are a nonlinear function of the power injections, as well as the network parameters, when determining the regulation power provided by each DER. 
As such, in the ODCP formulation, we approximate the losses using the so-called loss factors (LFs)---linear sensitivities of the total system losses with respect to changes in power injections---that explicitly capture the impacts of both active and reactive power injections on system losses.
LFs are conventionally computed from power flow models and have been applied in locational marginal price computation in electricity markets \cite{LF1} and sizing and allocation of DERs in power distribution systems \cite{LF2}.
In the proposed framework, the LF will be estimated in an online fashion using measurements acquired from the system, and updated in real-time so as to adapt to various phenomena that impact the operation of the power system, e.g., changes in system operating point.


\section{Preliminaries} \label{sec:prelim}

In this section, we first introduce the power distribution system model used throughout this paper.
Then, we describe the problem of coordinating a set of DERs to collectively provide frequency regulation services to a bulk power grid.

\subsection{Power Distribution System Model} \label{sec:model}

Consider a balanced power distribution system that consists of a set of buses indexed by $\tilde{\mathcal{N}} = \{0, 1, \cdots, N\}$.
Assume that bus $0$ corresponds to a substation bus, which is the only connection of the distribution system to a sub-transmission (or transmission) network. 
Further, assume that the bus $0$ is an ideal voltage source.

Without loss of generality, assume there is one DER and one load at every bus except for bus $0$, where there is no DER or load.
Let $P_i^g$ and $Q_i^g,~i=1,2,\dots,N,$ respectively denote the active and reactive power injections from the DER at bus $i$ (referred to as DER $i$).
Similarly, let $P_i^d$ and $Q_i^d$ respectively denote the  active and reactive power withdrawals of the load at bus $i$ (referred to as load $i$).
Typically, there will be some reactive power control resources in the power distribution system to maintain the bus voltages within a desired range \cite{volt_control}.
Assume the reactive power injections from such resources are included in the  $ {Q}^d_i$'s.
Let $P_i$ and $Q_i$ denote the respective net active and reactive power injections at bus~$i$, and $P^t$ and $Q^t$ the respective active and reactive power injections into the power distribution system from the bulk power grid.
Then, $P_0 = P^t$, $Q_0 = Q^t$, and  for any $i \in \mathcal{N}$, where $\mathcal{N} = \{1, \cdots, N\}$, $P_i = P_i^g - P_i^d$, $Q_i = Q_i^g - Q_i^d$.
Define $\bm{P} =\bm{P}^g-\bm{P}^d$, where $\bm{P}^g = [P_1^g, \cdots, P_N^g]\trans$, and $\bm{P}^d = [P_1^d, \cdots, P_N^d]\trans$. Similarly, define  $\bm{Q} =\bm{Q}^g -\bm{Q}^d $, where  $\bm{Q}^g = [Q_1^g, \cdots, Q_N^g]\trans$, and $\bm{Q}^d = [Q_1^d, \cdots, Q_N^d]\trans$. 
Then, we can conceptually express $P^t$ as a function of $\bm{P}$ and $\bm{Q}$ as
\begin{equation} \label{eq:h_function}
    P^t = h(\bm{P}, \bm{Q}).
\end{equation}

For given $\bm{P}$ and $\bm{Q}$, the total system losses, denoted by $\ell(\bm{P}, \bm{Q})$, are given by
\begin{equation} \label{eq:loss}
    \ell(\bm{P}, \bm{Q}) = h(\bm{P}, \bm{Q}) + \ones_N\trans \bm{P},
\end{equation}
where $\ones_N$ is an $N$-dimensional all-ones vector.
The partial derivatives of the total system losses with respective to active and reactive power injections at each bus are referred to as the loss factors (LFs).
Let $\Lambda_i^p$ and $\Lambda_i^q$ respectively denote the partial derivatives of the total system losses with  respect to the active and reactive net power injections at bus $i$, $i \in \mathcal{N}$; we refer to $\Lambda_i^p$ as the active LF at bus $i$, and $\Lambda_i^q$ as the reactive LF at bus $i$.
Define $\bm{\Lambda}^p = [\Lambda_1^p, \cdots, \Lambda_N^p]\trans$ and $\bm{\Lambda}^q = [\Lambda_1^q, \cdots, \Lambda_N^q]\trans$.
Then, it follows from \eqref{eq:loss} that 
\begin{equation} \label{eq:LF_def}
\begin{array}{lll}
    (\bm{\Lambda}^p)\trans \coloneqq  \left.\dfrac{\partial \ell}{\partial \bm{P}}\right\vert_{(\bm{P}, \bm{Q})} =  \left.\dfrac{\partial h}{\partial \bm{P}}\right\vert_{(\bm{P}, \bm{Q})} + \ones_N\trans, \\
    (\bm{\Lambda}^q)\trans \coloneqq \left.\dfrac{\partial \ell}{\partial \bm{Q}}\right\vert_{(\bm{P}, \bm{Q})} = \left.\dfrac{\partial h}{\partial \bm{Q}}\right\vert_{(\bm{P}, \bm{Q})},
\end{array}
\end{equation}
where the partial derivatives of $h$ with respect to $\bm{P}$ and $\bm{Q}$ are row vectors.
Expressions for $\bm{\Lambda}^p$ and $\bm{\Lambda}^q$ can be obtained using a power flow model; due to the space limitation, the derivation of model-based LFs are not presented here.

Therefore, for given small $\Delta \bm{P}$ and $\Delta \bm{Q}$ (with respect to some $\bm{P}$ and $\bm{Q}$), the incremental total system losses associated with $\Delta \bm{P}$ and $\Delta \bm{Q}$ can be approximated by
\begin{equation} \label{eq:inc-loss-LF}
    \Delta \ell(\bm{P}, \bm{Q}) \approx (\bm{\Lambda}^p)\trans \Delta \bm{P} + (\bm{\Lambda}^q)\trans \Delta \bm{Q}.
\end{equation}
Then, by using  \eqref{eq:h_function} and \eqref{eq:LF_def}, we can obtain the change in $P^t$, denoted by $\Delta P^t$, as follows:
\begin{equation} \label{eq:P1-P_1}
\begin{array}{lll}
	\Delta P^t & \approx & \left.\dfrac{\partial h}{\partial \bm{P}}\right\vert_{(\bm{P}, \bm{Q})} \Delta \bm{P} + \left.\dfrac{\partial h}{\partial \bm{Q}}\right\vert_{(\bm{P}, \bm{Q})} \Delta \bm{Q} \\
    & = & (\bm{\Lambda}^p - \ones_N)\trans \Delta \bm{P} + (\bm{\Lambda}^q)\trans \Delta \bm{Q}.
\end{array}
\end{equation}

It is clear from \eqref{eq:P1-P_1} that the impacts from reactive power injections on the total system losses cannot be neglected.
When active power injections change, the net reactive power change, $\Delta \bm{Q}$, may change accordingly based on some specific rules (possibly as a result of a feedback control action).
Both $\Delta \bm{P}$ and $\Delta \bm{Q}$ will lead to changes in system losses.
Thus, in order to determine the incremental total system losses after an active power injections change, it is necessary to know the reactive power control policies implemented throughout the system.
If $\Delta \bm{P}$ is small, we can approximately represent $\Delta \bm{Q}$ as $\Delta \bm{Q} \approx \bm{\Phi} \Delta \bm{P}$,
where $\bm{\Phi} \in \real^{N \times N}$ is the reactive power response sensitivity matrix at the operating point defined by $(\bm{\theta}, \bm{V}, \bm{P}, \bm{Q})$.
Using this approximation, together with \eqref{eq:inc-loss-LF} and \eqref{eq:P1-P_1}, $\Delta \ell(\bm{P})$ and $\Delta P^t $ can be represented as a linear function of $\Delta \bm{P}$ as follows:
\begin{align}
\Delta \ell(\bm{P}, \bm{Q}) &\approx \bm{\Lambda}\trans \Delta \bm{P}, \\ 
\Delta P^t &\approx (\bm{\Lambda} - \ones_N)\trans \Delta \bm{P}, \label{eq:loss-P-relation}
\end{align}
where $\bm{\Lambda} \coloneqq \bm{\Lambda}^p + \bm{\Phi}\trans \bm{\Lambda}^q$; we refer to the entries of $\bm{\Lambda}$ as the total loss factors.\footnote{It is important to emphasize that the total LFs depend strongly on the reactive power control policies.  In the simple case where no reactive power control is employed, i.e., $\bm{\Phi} = \bm{0}$, the total LFs are identical to the active LFs. However, in general, $\bm{\Phi}$ can be hardly derived from models.}

\subsection{DER Coordination for Frequency Regulation}

Let $\bm{P}^{g0}_i$ ($\bm{P}^{d0}_i$) denote the nominal power injection (withdrawal) from DER (load) $i$, $i \in \mathcal{N}$,  and define $\bm{P}^{g0}=[P_1^{g0}, \cdots, P_N^{g0}]\trans$ ($\bm{P}^{d0}=[P_1^{d0}, \cdots, P_N^{d0}]\trans$). 
Let $P^{t0}$ denote the nominal active power injection from the bulk power grid.  

Let $r(t)$ denote the value of some regulation signal, $r$, sent by the bulk system operator to the distribution system at time~$t$. Usually, the operator sends a single value every $\Delta T$ units of time,  e.g.,   $\Delta T=2$~s. Thus, in reality, the signal  $r(t)$ is  piecewise constant, taking  some  value $r[k] \in \mathbb{R}$ for $t$ in the interval $(k \Delta T, (k+1) \Delta T],~k=0,1,\dots$,  referred to as time interval~$k$.
Let $P_i^{g0}[k]$ denote the nominal power injected by DER~$i$ in time interval $k$, and let $P_i^{g}[k]$ denote the actual power injected by DER~$i$ in the same time interval. 
Then, we can define the regulation power provided by DER $i$ in time interval $k$, as $p_i^g[k]:=P_i^g[k]-P_i^{g0}[k]$,\footnote{The nominal power injection may vary with time; hence its dependence of $k$. For example, $P_i^{g0}[k]$ may increase/decrease during a period so as to follow forecasted  load changes.} and  the corresponding regulation power vector for the time  interval~$k$ as $\bm{p}^g[k] = [p_1^g[k], \cdots, p_N^g[k]]\trans$.
Let $\overline{p}_i^g $ and $\underline{p}_i^g $ respectively denote the maximum up and down regulation capacities of DER $i$, and define $\overline{\bm{p}}^g = [\overline{p}_1^g, \cdots, \overline{p}_N^g]\trans$ and $\underline{\bm{p}}^g = [\underline{p}_1^g, \cdots, \underline{p}_N^g]\trans$.
Then, at time instant $k\Delta T,~k=0,1,\dots$,  referred to as time instant $k$, the problem is to determine $\bm{p}^g[k]$ satisfying the following two constraints: 
\begin{itemize}
\item[\textbf{[C1.]}] $ \underline{\bm{p}}^g \leq \bm{p}^g[k] \leq \overline{\bm{p}}^g$, i.e., DER up and down regulation capacity constraints are satisfied; and 
\item[\textbf{[C2.]}] $P^t[k] = P^{t0}[k] - r[k]$, i.e., the total amount of regulation power provided to the bulk grid is $r[k]$,\footnote{We adopt the convention that if $r[k]$ is positive, then the distribution system must regulate upwards, i.e., the total amount of power injected into the distribution system must be decreased.}
\end{itemize}
while all loads within the power distribution system are balanced, i.e., equation \eqref{eq:loss-P-relation}, which is used to approximately model the power balance in the system, is satisfied for some~$\bm{Q}$.
Note that $r[k]$ and all the measurements are obtained at time instant $k$, and the loads are assumed to be constant for the duration of the interval $(k\Delta T, (k+1) \Delta T]$. Also, the computation of $\bm{p}^g[k]$ is assumed to be instantaneous.
The DERs will be instructed to provide $\bm{p}^g[k]$ immediately after the computation is completed, and their power outputs will remain fixed for the duration of the interval $(k\Delta T, (k+1) \Delta T]$.

In a lossless power distribution system, the sum of $p_i^g[k]$'s must equal to $r[k]$; however, this is not the case in a lossy power distribution system.
As such, it is important to explicitly take into account the system losses in order for the DERs to provide the exact amount of regulation power to the bulk power grid. Also, while there may exist multiple ways to choose a   $\bm{p}^g[k]$ that solves the DER coordination problem, we can select one that is optimal with respect to some objective function.
For example, we may select a $\bm{p}^g[k]$ that minimizes the incremental total system losses or costs as in \cite{DER_coord}, or one that minimizes the norm of  the regulation power vector, or a combination of these various objectives.

\section{Adaptive DER Coordination Framework} \label{sec:DER_coord}

In this section, we propose an adaptive coordination framework for DERs to provide frequency regulation.
We first propose a measurement-based method to estimate the LFs, which is then used to develop a formulation of the ODCP. 

\subsection{Recursive Loss Factor Estimator} \label{sec:RLFE}

We next propose a measurement-based approach for estimating LFs, which is built upon the recursive weighted least-squares (RWLS) estimation method. 
Suppose at time instant $k$, we have $k$  consecutive measurements of $P^t$ and $\bm{P}$, denoted by $P^t[0], \cdots, P^t[k-1]$, $\bm{P}[0], \cdots, \bm{P}[k-1]$, respectively (measurements $P^t[l]$ and $\bm{P}[l]$, $l = 0, 1, \cdots, k-1$, are obtained at time instant $l$). 
Define the change in the power injection vector at time instant $k$ as $\Delta \bm{P}[l] = \bm{P}[l] - \bm{P}[l-1]$, where 
\begin{align*}
	\bm{P}[l] & =   \bm{P}^{g0}[l] + \bm{p}^g[l] - \bm{P}^{d}[l], \\
	\bm{P}[l-1] & =   \bm{P}^{g0}[l-1] + \bm{p}^g[l-1] - \bm{P}^{d}[l-1].
\end{align*}
Let  $ \Delta P^t[l] $ denote the change in the active power injection into the network from the bulk power grid at time instant $l$ that results from $\Delta \bm{P}[l]$, i.e., $\Delta P^t[l] = P^t[l] - P^t[l-1]$. Then, for sufficiently small $\Delta \bm{P}[l]$, and assuming $\bm{\Lambda}$ remains relatively constant over time if there are no topology changes,  it follows from \eqref{eq:loss-P-relation} that
\begin{equation} \label{eq:lse_base}
	\Delta \bm{P}\trans[l] (\bm{\Lambda} - \ones_N) \approx \Delta P^t[l]. 
\end{equation}
Then, by stacking the equation in \eqref{eq:lse_base} for $l=0,1,\cdots,k-1$, we obtain   the following system of equations:
\begin{equation} \label{eq:lse_LF}
    \begin{bmatrix}
    \Delta \bm{P}\trans[1] \\ 
    \vdots \\ 
    \Delta \bm{P}\trans[k-1]
    \end{bmatrix}
    (\bm{\Lambda} - \ones_N) \approx 
    \begin{bmatrix}
    \Delta P^t[1] \\ \vdots \\ \Delta P^t[k-1]
    \end{bmatrix}.
\end{equation}
When $k$ is sufficiently large, \eqref{eq:lse_LF} becomes   overdetermined, and thus, we can use weighted least-squares (WLS)   to  obtain an estimate of $\bm{\Lambda}$ at instant $k-1$, denoted by $\hat{\bm{\Lambda}}[k-1]$,     as follows:
\begin{equation}
    \hat{\bm{\Lambda}}[k-1] = \ones_N + \bm{R}[k-1] \bm{A}\trans[k-1] \bm{W}[k-1] \bm{b}[k-1],
\end{equation}
where  $\bm{W}[k-1] \in \mathbb{R}^{(k-1)\times (k-1)}$ is a positive definite symmetric weight matrix, and 
\[
\bm{A}[k-1] = 
\begin{bmatrix}
\Delta \bm{P}\trans[1] \\ 
\vdots \\ 
\Delta \bm{P}\trans[k-1]
\end{bmatrix},
\bm{b}[k-1] = \begin{bmatrix} \Delta P^t[1] \\ \vdots \\ \Delta P^t[k-1] \end{bmatrix},
\]
\[
\bm{R}[k-1] = \left(\bm{A}\trans[k-1] \bm{W}[k-1] \bm{A}[k-1] \right)^{-1}.
\]

Let $\bm{W}[k-1] = \text{diag}\{\gamma^0, \gamma^1, \cdots, \gamma^{k-2}\}$ where $\gamma \in (0, 1]$ is the forgetting factor for the measurements.
When a new set of measurements, $P^t[k]$ and $\bm{P}[k]$, becomes available at $k$, we can update the LF estimate from $  \hat{\bm{\Lambda}}[k-1]$ to $  \hat{\bm{\Lambda}}[k]$  as follows \cite{rls}: 
\begin{equation} \label{eq:rwls}
    \begin{array}{lll}
        \hat{\bm{\Lambda}}[k] & = & \hat{\bm{\Lambda}}[k-1] + \bm{R}[k] \Delta \bm{P}[k](\Delta P^t[k] \\
        & & - \Delta \bm{P}\trans[k] (\hat{\bm{\Lambda}}[k-1] - \ones_N) ),
    \end{array}
\end{equation}
where
\begin{equation*}
    \bm{R}[k] = \dfrac{1}{\gamma} \left (\bm{R}[k-1] - \dfrac{\bm{R}[k-1] \Delta \bm{P}[k] \Delta \bm{P}\trans[k] \bm{R}[k-1]}{\gamma + \Delta \bm{P}\trans[k] \bm{R}[k-1] \Delta \bm{P}[k]} \right). 
\end{equation*}
Thus, using \eqref{eq:rwls}, the LFs can be updated using new measurements with little computational effort.

\subsection{LF-based ODCP Formulation}

Using the estimated LFs,  we can now develop a formulation for the ODCP.
At time instant $k$, the total LFs estimated at time instant $k-1$, $\hat{\bm{\Lambda}}[k-1]$, is available.
Then, for sufficiently small $\Delta \bm{P}[k]$, it follows from  \eqref{eq:lse_base} that
\begin{align} \label{eq:incremental_model}
\Delta P^t[k] \approx (\hat{\bm{\Lambda}}[k-1] - \ones_N)\trans \Delta \bm{P}[k].
\end{align}

Now, given $r[k]$, we  would like to choose $ \bm{p}^g[k] $ in some optimal fashion so that that $ P^t[k] $ tracks the regulation signal, i.e., $ P^t[k]=P^{t0}[k]-r[k]$. Thus, by using the incremental model in \eqref{eq:incremental_model}, we have that 
\begin{align}
\Delta P^t[k]&=P^t[k]- P^t[k-1] \nonumber \\
&= (P^{t0}[k] - r[k]) - P^t[k-1] \nonumber \\
& =(\hat{\bm{\Lambda}}[k-1] - \ones_N)\trans \Delta \bm{P}[k],
\end{align}
with
\begin{align}
\Delta \bm{P}[k] & = \bm{P}[k] - \bm{P}[k-1] \nonumber \\
    & = (\bm{P}^{g0}[k] + \bm{p}^g[k] - \bm{P}^{d}[k]) \nonumber \\
    & - (\bm{P}^{g0}[k-1] + \bm{p}^g[k-1] - \bm{P}^d[k-1]).
\end{align}

Then, by choosing  the minimization of a weighted sum of the incremental losses and the $L_2$-norm (denoted by $\norm{\cdot}$) of the regulation power vector as the optimality criterion,  the ODCP to be solved at time instant $k$ can be formulated as follows:
\begin{equation*} 
\bm{p}^g[k] = \argmin_{\bm{z}} ~ \hat{\bm{\Lambda}}\trans[k-1] \Delta \bm{P}[k] + \frac{\rho}{2} \norm{\bm{z}}^2
\end{equation*}
subject to 
\begin{subequations} \label{eq:ODCP}
\begin{equation} \label{eq:c1}
    (P^{t0}[k] - r[k]) - P^t[k-1] = (\hat{\bm{\Lambda}}[k-1] - \ones_N)\trans \Delta \bm{P}[k],
\end{equation}
\begin{equation} \label{eq:c2}
\begin{array}{lll}
    \Delta \bm{P}[k] & = & (\bm{P}^{g0}[k] + \bm{z} - \bm{P}^{d}[k]) \\
    & & - (\bm{P}^{g0}[k-1] + \bm{p}^g[k-1] - \bm{P}^d[k-1]),
\end{array}
\end{equation}
\begin{equation} \label{eq:c3}
    \underline{\bm{p}}^g \leq \bm{z} \leq \overline{\bm{p}}^g.
\end{equation}
\end{subequations}
where $\rho \geq 0$, $\bm{z} = [z_1, \cdots, z_N]\trans$ with $z_i$ being the regulation power provided by DER $i$, which is to be determined.
Note that at time instant $k$, $P^t[k-1]$, $\bm{P}^d[k-1]$, $\bm{P}^{g0}[k-1]$, $P^{t0}[k]$, $\bm{P}^d[k]$, $\bm{P}^{g0}[k]$, and $r[k]$ are known quantities.
The ODCP is a quadratic program with one equality constraint, which can be solved efficiently.

\subsection{Interaction between LF Estimator and ODCP Solver}
At time instant $k$, the LF estimator provides $\hat{\bm{\Lambda}}[k-1]$, i.e., the estimate of the total LFs based on some measurements taken up to time instant $k-1$.
Then, the vector of estimated LFs, $\hat{\bm{\Lambda}}[k-1]$, is sent to the ODCP solver.
The ODCP solver also receives loads $\bm{P}^d[k]$, the requested regulation power $r[k]$, and solves the ODCP to determine $\bm{p}^g[k]$.
Then, DERs are instructed to change their set-points so as to provide $\bm{p}^g[k]$ for $k\Delta T  < t \leq   (k+1)\Delta T$.
A new set of measurements will be available after the DER set-points are modified.
The new measurements will be used in the  LF estimator to dynamically update the estimated values of the total LFs and obtain $\hat{\bm{\Lambda}}[k]$, which will be used in ODCP to compute the optimal regulation power from DERs for the next time interval.


\section{Numerical Simulation} \label{sec:simu}

In this section, we illustrate the application of the proposed framework through numerical simulations.
The accuracy of the estimated LFs, as well as the performance of the framework, is studied for a case where the nominal loads are kept constant, and a case where there are load changes.

\subsection{Simulation Setup}

A modified 33-bus power distribution system from \cite{33bus} is used for all numerical simulations. 
There are three DERs, located at buses 12, 25, and 33, respectively, with their respective capacities being $2300$ kW, $1500$ kW, and $1200$ kW, and their regulation capacity equal to 10\% of its total capacity.
We assume the reactive power control at buses 1 and 12 aims to maintain a constant voltage magnitude of $1$ p.u. while no reactive  power control is employed at any other buses.
Throughout the simulation, we set $\rho = 1$ and $\gamma = 0.97$.

The power demanded by load $i$ is generated using $P_i^d[k] = P_i^{d0} (1 + \nu_i)$, where $P_i^{d0}$ denotes the demand nominal value, and $\nu_i$ is a zero-mean Gaussian random variable with standard deviation being $\sigma_i = 0.01$. 
The regulation signal is taken from PJM \cite{PJM_AGC_signal} and is updated every 2 seconds; correspondingly, we update the DER set-points every 2 seconds.

\begin{figure}[!t]
\centering
\includegraphics[width=3.5in]{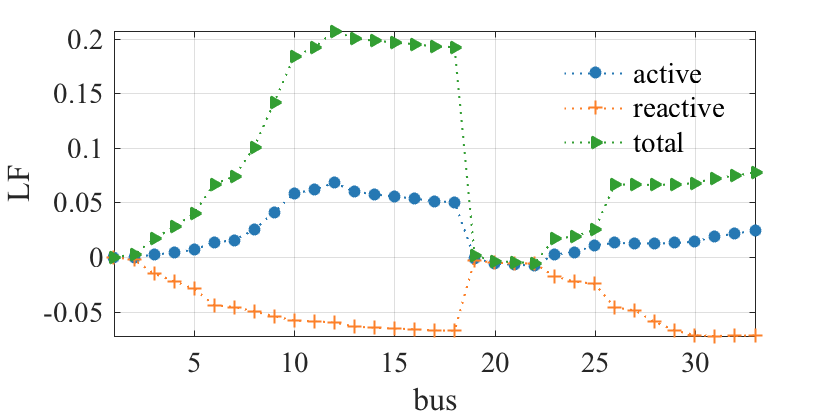}
\caption{Actual LFs at nominal loads.}
\label{fig:true_LFs}
\end{figure}

The actual LFs are obtained by manually perturbing the system and their values at the nominal operating point are shown in Fig. \ref{fig:true_LFs}.
As can be seen from Fig. \ref{fig:true_LFs}, the active LFs, reactive LFs, and the total LFs have distinct values.
Therefore, absent the knowledge of $\bm{\Phi}$, we can hardly obtain the accurate total LFs from the model.

\subsection{LF Estimation Accuracy} \label{sec:LF_est}

\begin{figure}[!t]
\centering
\includegraphics[width=3.5in]{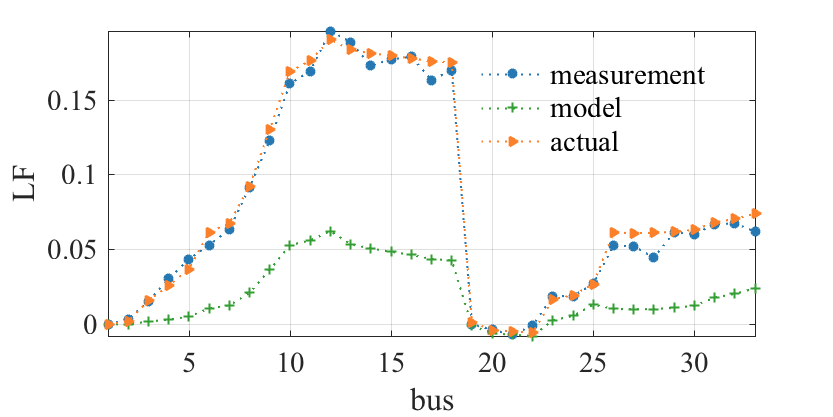}
\caption{LF estimation results at 0s.}
\label{fig:LF_est}
\end{figure}

\begin{figure}[!t]
\centering
\includegraphics[width=3.5in]{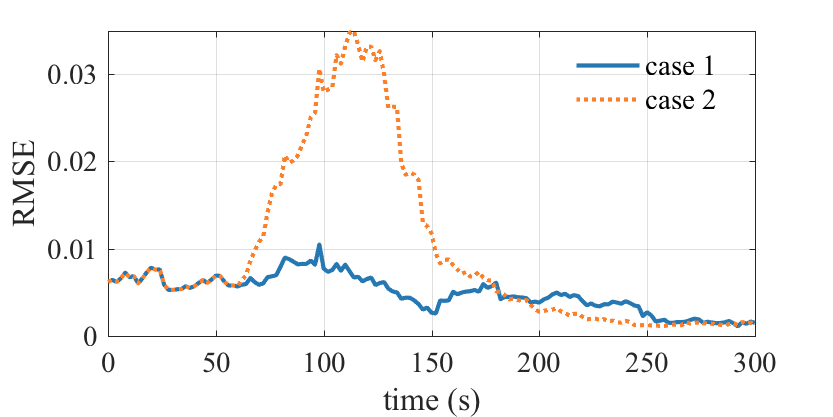}
\caption{LF estimation errors.}
\label{fig:LF_error_compare}
\end{figure}

The accuracy of the LF estimation is measured by the root mean squared error (RMSE), denoted by $\varepsilon$.
Assume we have an initial estimation of $\hat{\bm{\Lambda}}$ and $\bm{R}$ at $t=0$~s obtained by solving \eqref{eq:lse_LF}.
The estimation results are shown in Fig. \ref{fig:LF_est}.
In the model-based approach, the active LFs are used as estimates for the total LFs.
At $t=0$~s, $\varepsilon[0] = 0.0062$ in the measurement-based approach, which is relatively small.
Yet, we have $\varepsilon[0] = 0.0758$ for the model-based approach, which is one order of magnitude greater than the values obtained using the measurement-based approach.
This is expected since the impacts of reactive LFs on the total LFs are ignored in the model-based approach due to the lack of knowledge on $\bm{\Phi}$.
After the first estimation, the RWLS algorithm is used to dynamically update $\hat{\bm{\Lambda}}$ and $\bm{R}$.
The RMSEs during $t=0$~s to $t=300$~s are shown as the solid line (case 1) in Fig. \ref{fig:LF_error_compare}.
The average RMSE in the measurement-based approach is $0.0049$.

To illustrate the adaptivity of the measurement-based approach, we simulate a case where the nominal loads increase linearly by $20$\% from $t=60$~s to $t=120$~s.
Under this setup, the RMSEs during the time interval $(0,300]$~s are presented as the dashed line (case~2) in Fig. \ref{fig:LF_error_compare}.
The average RMSE in the measurement-based approach is $0.0096$, while the maximum RMSE is $0.0349$.
The RMSE starts to increase when the operating point begins to change at $t=60$~s and starts to decrease when the changes end at $t=120$~s.
This is intuitively reasonable since the during the time when the operating point changes, the majority of existing measurements provide information on old operating points, based on which we can hardly get accurate estimates at the current operating point.
After the operating point change ends, the new measurements provide more information on the new operating point and the impacts from the measurements on the old one decay.

\subsection{Quantification of Frequency Regulation Performance}

The proposed framework is compared with the participation factor (PF) based coordination approach, which are non-negative real numbers that sum up to 1 and are proportional to the regulation capacity of each DER.
The total incremental changes in the active loads and the requested regulation power are allocated to each DER based on the PFs.
The frequency regulation performance is measured by a score denoted by $S$.
The value of $S$ for the interval $(k\Delta T, (k+1)\Delta T]$ is computed as $S[k] = 1 - \frac{\sum_{l=0}^k |r_m[l] - r[l]|}{\sum_{l=0}^k |r[l]|}$, where $r_m[l]=P^{t0}[l]-P^t[l]$ is the actual regulation power at time instant $l$.

In the base case, the average performance score obtained using the PF-based approach is $0.8756$; that obtained using the LF-based approach with estimated total LFs is $0.9991$; and that obtained using the LF-based approach with actual total LFs is $0.9992$.
Obviously, the LF-based approaches perform much better than the PF-based approach where system losses are ignored.
Moreover, the estimated-LF based approach perform almost equally good as the actual-LF based approach.

In the case with same operating point changes described in Section \ref{sec:LF_est}, the average performance score using LF-based approach with estimated total LFs is $0.9988$, and that obtained using the LF-based approach with actual total LFs is $0.9992$.
The slight decrease in the performance score obtained using the estimated LFs are caused by the decrease in the estimation accuracy of the total LFs.
Yet, the performance score decrease is negligible, showing the adaptivity of the propose framework.

\section{Concluding Remarks} \label{sec:con}

In this paper, we proposed an adaptive coordination framework for DERs in a lossy power distribution system to collectively provide frequency regulation services to a bulk power grid.
The framework consists of a measurement-based recursive LF estimator to capture the impacts of both active and reactive power injections on system losses, and an ODCP solver that determines the optimal DER power outputs using the estimated LFs.
The inherent nature of the estimator makes it adaptive to system condition changes.
Numerical simulation validated the effectiveness of the proposed framework in coordinating the DERs for frequency regulation.

\bibliographystyle{IEEEtran}
\bibliography{ref}

\begin{thebibliography}{10}
\providecommand{\url}[1]{#1}
\csname url@samestyle\endcsname
\providecommand{\newblock}{\relax}
\providecommand{\bibinfo}[2]{#2}
\providecommand{\BIBentrySTDinterwordspacing}{\spaceskip=0pt\relax}
\providecommand{\BIBentryALTinterwordstretchfactor}{4}
\providecommand{\BIBentryALTinterwordspacing}{\spaceskip=\fontdimen2\font plus
\BIBentryALTinterwordstretchfactor\fontdimen3\font minus
  \fontdimen4\font\relax}
\providecommand{\BIBforeignlanguage}[2]{{%
\expandafter\ifx\csname l@#1\endcsname\relax
\typeout{** WARNING: IEEEtran.bst: No hyphenation pattern has been}%
\typeout{** loaded for the language `#1'. Using the pattern for}%
\typeout{** the default language instead.}%
\else
\language=\csname l@#1\endcsname
\fi
#2}}
\providecommand{\BIBdecl}{\relax}
\BIBdecl

\bibitem{wind_impacts}
Y.~V. Makarov, C.~Loutan, J.~Ma, and P.~de~Mello, ``Operational impacts of wind
  generation on california power systems,'' \emph{IEEE Trans. Power Syst.},
  vol.~24, no.~2, pp. 1039--1050, May 2009.

\bibitem{ES}
F.~Zhang, Z.~Hu, X.~Xie, J.~Zhang, and Y.~Song, ``Assessment of the
  effectiveness of energy storage resources in the frequency regulation of a
  single-area power system,'' \emph{IEEE Trans. on Power Syst.}, vol.~32,
  no.~5, pp. 3373--3380, Sept. 2017.

\bibitem{performance}
B.~Xu, Y.~Dvorkin, D.~S. Kirschen, C.~A. Silva-Monroy, and J.~P. Watson, ``A
  comparison of policies on the participation of storage in u.s. frequency
  regulation markets,'' in \emph{Proc. of IEEE Power and Energy Society General
  Meeting}, July 2016, pp. 1--5.

\bibitem{DER_coord2}
A.~D. Dom\'{i}nguez-Garc\'{i}a and C.~N. Hadjicostis, ``Coordination and
  control of distributed energy resources for provision of ancillary
  services,'' in \emph{Proc. of IEEE International Conference on Smart Grid
  Communications}, Oct. 2010, pp. 537--542.

\bibitem{DER_model1}
S.~S. Guggilam, C.~Zhao, E.~Dall'Anese, Y.~C. Chen, and S.~V. Dhople, ``Primary
  frequency response with aggregated ders,'' in \emph{Proc. of American Control
  Conference}, May 2017, pp. 3386--3393.

\bibitem{meas_est}
Y.~C. Chen, A.~D. Dom\'inguez-Garc\'ia, and P.~W. Sauer, ``Measurement-based
  estimation of linear sensitivity distribution factors and applications,''
  \emph{IEEE Trans. Power Syst.}, vol.~29, no.~3, pp. 1372--1382, May 2014.

\bibitem{meas_est2}
J.~Zhang, X.~Zheng, Z.~Wang, L.~Guan, and C.~Y. Chung, ``Power system
  sensitivity identification -- inherent system properties and data quality,''
  \emph{IEEE Trans. Power Syst.}, vol.~PP, no.~99, pp. 1--1, 2016.

\bibitem{LF1}
E.~Litvinov, T.~Zheng, G.~Rosenwald, and P.~Shamsollahi, ``Marginal loss
  modeling in lmp calculation,'' \emph{IEEE Trans. Power Syst.}, vol.~19,
  no.~2, pp. 880--888, May 2004.

\bibitem{LF2}
N.~Acharya, P.~Mahat, and N.~Mithulananthan, ``An analytical approach for dg
  allocation in primary distribution network,'' \emph{Int. J. Elec. Power},
  vol.~28, no.~10, pp. 669--678, Dec. 2006.

\bibitem{volt_control}
B.~A. Robbins, C.~N. Hadjicostis, and A.~D. Dom\'{i}nguez-Garc\'{i}a, ``A
  two-stage distributed architecture for voltage control in power distribution
  systems,'' \emph{IEEE Trans. Power Syst.}, vol.~28, no.~2, pp. 1470--1482,
  May 2013.

\bibitem{DER_coord}
H.~Xu, S.~C. Utomi, A.~D. Dom\'{i}nguez-Garc\'{i}a, and P.~W. Sauer,
  ``Coordination of distributed energy resources in lossy networks for
  providing frequency regulation,'' in \emph{Proc. of X Bulk Power Systems
  Dynamics and Control Symposium}, Aug. 2017.

\bibitem{rls}
S.~S. Haykin, \emph{Adaptive Filter Theory}.\hskip 1em plus 0.5em minus
  0.4em\relax Prentice Hall, 1996.

\bibitem{33bus}
M.~E. Baran and F.~F. Wu, ``Network reconfiguration in distribution systems for
  loss reduction and load balancing,'' \emph{IEEE Trans. Power Del.}, vol.~4,
  no.~2, pp. 1401--1407, Apr. 1989.

\bibitem{PJM_AGC_signal}
\BIBentryALTinterwordspacing
PJM. (2017) {RTO} regulation signal data. [Online]. Available:
  \url{http://www.pjm.com/markets-and-operations/ancillary-services.aspx}
\BIBentrySTDinterwordspacing

\end{thebibliography}

\end{document}